\documentclass[12pt, letterpaper]{article}
\usepackage[utf8]{inputenc}
\usepackage{amsmath}
\usepackage{graphicx}
\usepackage{authblk}
\usepackage{amsthm}
\graphicspath{ {./images/} }
\usepackage[utf8]{inputenc}
\usepackage{amssymb}
\usepackage{multicol}
\usepackage{dcolumn}
\usepackage{changepage}

\newtheorem{lemma}{Lemma}
\newtheorem{definition}{Definition}
\newtheorem{corollary}{Corollary}
\newtheorem{theorem}{Theorem}

\begin{document}


\title{A note on asymptotic cones of graph-adapted smocked spaces}

\author[1,2,3]{Hollis Williams\thanks{holliswilliams@hotmail.co.uk}}

\affil[1]{Department of Physics, University of Warwick, Coventry CV4 7AL, United Kingdom}
\affil[2]{Theoretical Sciences Visiting Program, Okinawa Institute of Science and Technology Graduate University, Onna-son, Okinawa 904-0495, Japan}
\affil[3]{\textit{Current address:} School of Engineering, Westlake University, Hangzhou, Zhejiang 310030, China}

\maketitle

\begin{abstract}
Smocked spaces, introduced by Sormani and collaborators as a generalization of pulled thread spaces, provide a broad class of metric quotients of Euclidean space. In this note we investigate their large-scale geometry via periodic graph models.  We introduce the notion of a graph-adapted smocked realization of a periodic weighted graph and establish a uniform additive distortion estimate between the smocked metric and the underlying graph metric. As a consequence, graph-adapted smocking preserves stable norms and asymptotic cones. Combining this with classical homogenization results for periodic graphs, we show that the tangent cone at infinity of a graph-adapted smocked space is determined by the stable norm of the associated periodic graph, which implies that every centrally symmetric rational polyhedral norm arises as the unique tangent cone at infinity of a graph-adapted smocked space. This establishes a connection between stable norm theory and the asymptotic geometry of smocked spaces.


\end{abstract}

\textbf{Keywords}: Smocked spaces, Gromov-Hausdorff convergence, metric geometry, tangent cones.
\newline
\newline
\noindent
53C23, 54E35; Secondary 28A33, 51F30.

\section{Introduction}
\label{}

\noindent
Smocked spaces were introduced by Sormani et al. as a generalization of pulled thread spaces [1, 2]. Roughly speaking, one begins with Euclidean space and collapses a collection of pairwise disjoint compact connected subsets (known as stitches) to points. The resulting quotient metric spaces exhibit a rich geodesic structure and provide a flexible class of examples in metric geometry.  One of the main questions raised in the original work on these spaces concerns their asymptotic geometry. In particular, Sormani et al. asked whether every finite-dimensional normed vector space can arise as the tangent cone at infinity of a smocked space.  This question has some bearing on the topic of Riemannian manifolds with lower scalar curvature bounds, since when studying curvature under convergence, the limiting spaces are not always smooth manifolds and it is important to understand the metric structures of the large-scale limits [3].  Tangent cones at infinity provide one method to extract this large-scale geometry and although a number of examples of these cones have been constructed and studied in detail for smocked spaces, a general realization theorem remains open.

The purpose of the present note is to investigate this question from the viewpoint of periodic graphs and stable norms. Stable norms associated to periodic graphs have been studied extensively in geometric group theory, homogenization theory and the geometry of periodic metrics. These norms describe the large-scale geometry of the graph and determine its asymptotic cone.  The main result which we prove is the following:

\begin{theorem}
\label{thm:main}

Let $G$ be a connected $\mathbb Z^n$--periodic weighted graph with stable norm $\|\cdot\|_{\mathrm{st}}$ on $\mathbb R^n$.  Let $X$ be a graph-adapted smocked realization of $G$
satisfying assumptions {\rm (A1)--(A3)}.  Then there exists a constant $C>0$ such that
\[
|d_X(x,y)-d_G(x,y)|\le C
\]
for all $x,y$ in the vertex set (or orbit space) of the realization.

Consequently:

\begin{enumerate}
\item $X$ and $G$ have the same stable norm,
\[
\|\cdot\|_{\mathrm{st}}^X
=
\|\cdot\|_{\mathrm{st}}^G;
\]

\item the asymptotic cones of $X$ and $G$ coincide;

\item if
\[
(G,\lambda^{-1}d_G,o)
\xrightarrow{\mathrm{pGH}}
(\mathbb R^n,F),
\]
then
\[
(X,\lambda^{-1}d_X,o)
\xrightarrow{\mathrm{pGH}}
(\mathbb R^n,F).
\]
\end{enumerate}

In particular, $X$ and $G$ have the same unique tangent cone at
infinity.
\end{theorem}

\noindent
This theorem states that an appropriately constructed smocking of a periodic graph $G$ is always a uniformly bounded additive distance from the underlying graph, implying that the two spaces have the same asymptotic geometry.  Combining this observation with results for stable norms of periodic graphs, we obtain a large class of normed spaces which can be realized as the tangent cone at infinity of suitably constructed smocked spaces:


 \begin{corollary}
\label{cor:polyhedral}
 Every centrally symmetric rational polyhedral norm on $\mathbb{R}^n$ arises as the unique tangent cone at infinity of a
graph-adapted smocked space.

\end{corollary}



 
 \noindent
Although this class of norms is quite restrictive, they can be used to approximate every norm on $\mathbb{R}^n$ arbitrarily well, so intuitively one may think of the result as saying that any finite-dimensional normed space can be well-approximated as the tangent cone of an appropriately constructed smocked space.   

The step to rational polyhedral norms uses a standard consequence of classical realization theory for stable norms (stated in the text as Lemma 1) that every centrally symmetric rational polyhedral norm on $\mathbb{R}^n$ can be realized as the stable norm of a connected $\mathbb{Z}^n$-periodic weighted graph (see Babenko-Balacheff and Jotz for closely related results) [4, 5].  Stable norm realization allows one to encode a broad and geometrically natural class of normed spaces by periodic weighted graphs and Theorem 1 then transfers the large-scale geometry of such a graph to a graph-adapted smocked space.  This opens up an interesting web of connections between stable norms on graphs and periodic geometric models realized by smocked spaces which could be explored further.



\section{Preliminaries}

\noindent
We will begin with some necessary definitions concerning the Gromov-Hausdorff distance.

\begin{definition}: Given two metric spaces $X$ and $Y$, the Gromov-Hausdorff distance is defined to be 

\[ d_{GH}(X,Y) = \text{inf} \bigg\{ d_H^Z(f(X),f(Y)) \bigg\} \]

\noindent
taken over all metric spaces $Z$ and isometric embeddings $f: X \rightarrow Z$ and $g: Y \rightarrow Z$, where $d_H^Z$ is the Hausdorff distance between subsets of $Z$:

\[ d_H^Z (A,B) =  \text{inf} \{ \epsilon > 0 : B \subset T_{\epsilon}(A) \: \text{and} \: A \subset T_{\epsilon}(B) \} ,\]

\noindent
and $T_{\epsilon}(A) = \{ x \in Z : d_Z(x,A) < \epsilon \}$. \end{definition}

\noindent
This distance has an associated notion of convergence known as (pointed) Gromov-Hausdorff convergence.

\begin{definition}: If one has a sequence of complete noncompact metric spaces $(X_j, d_j)$ and points $x_j \in X_j$, pointed Gromov-Hausdorff convergence is defined as follows:

\[ (X_j, d_j, x_j) \xrightarrow{pGH} (X_{\infty}, d_{\infty}, x_{\infty})    \]

\noindent
if and only if for every $R > 0$ the closures of the balls of radius $R$ in $X_j$ converge in the Gromov-Hausdorff sense as metric spaces with the restriction of the distance function to closures of balls in $X_{\infty}$:

\[ (\overline{B}_R (x_j), d_j) \xrightarrow{GH} (\overline{B}_R (x_{\infty}), d_{\infty}).  \]
\end{definition}

\noindent
A sequence of metric spaces can be obtained by taking a space and repeatedly rescaling it.
Although some spaces like the taxi space are isometric to their rescalings, this
is generally not the case. If a sequence or subsequence obtained by rescaling
converges in the Gromov-Hausdorff sense, the limit space obtained is known
as the tangent cone at infinity.

\begin{definition} A pointed metric space $(Y,d_Y,y_0)$ is a tangent cone at infinity
of a proper metric space $(X,d_X, x_0)$ if there exists a sequence
$\lambda_j\to\infty$ such that

\[
(X,\lambda_j^{-1}d_X,x_0)
\xrightarrow{\mathrm{pGH}}
(Y,d_Y,y_0).
\]
\end{definition}

We next provide the definition for a smocked space (the reader may think of such a space intuitively as the generalisation of a pulled thread space).  Roughly speaking, a smocked space $X$ is obtained by starting with a countable collection of pairwise disjoint compact connected subsets $\{I_j\}_{j\in J}$ in $\mathbb E^N$ and collapsing each $I_j$ to a point.  The induced length metric on the quotient defined by the quotient map $\pi:\mathbb E^N\to X$ is called the smocking metric.  For completeness, we will provide the formal definition below:

\begin{definition} A smocked space is a metric space $(X,d)$ defined as 

\[ X = \{ x: x \in \mathbb{E}^n \: \backslash \: S \} \cup I ,\]

\noindent
where $S$ is a smocking set defined by 

\[ S = \bigcup_{j \in J} I_j  \]

\noindent
for some index set $J$ and $I$ is a collection of countably many disjoint connected compact sets in $\mathbb{E}^N$:

\[I = \{ I_j : j \in J \} .\]

\noindent
Furthermore, define the smocking map $\pi: \mathbb{E}^N \rightarrow X$ by

\[ \pi(x) = \left\{
	\begin{array}{ll}
		x  & \mbox{if } x \in \mathbb{E}^n \: \backslash \: S \\
	 I_j  & \mbox{if } x \in I_j, \: j \in J\\
	\end{array}
\right. \]

\noindent 
The smocked distance function $d: X \times X \rightarrow [0, \infty)$ is defined as follows for $x,y \notin \pi(S)$ and stitches $I_m$ and $I_n$:

\[ d(x,y)   = \text{min} \{ d_0 (x,y), d_1 (x,y), ... \}  , \]

\[ d(x, I_n) = \text{min} \{ d_0 (x,z), d_1 (x,z), ... : z \in I_n \}     ,\]

\[ d( I_m, I_n) =  \text{min} \{ d_0 (z',z), d_1 (z',z), ... : z' \in I_m, z \in I_n \} ,\]

\noindent
where $d_k$ is the sum of lengths of segments which go between $k $ stitches: 

\[ d_k (v,w) = \text{min} \{ |v- z_1 | + \sum_{i=1}^{k-1} | z'_i - z_{i+1} | + |z'_k - w | : z_i, z'_i \in I_{j_i}, j_1 \neq ... \neq j_k \in J \}  ,\]

\noindent 
for any $v, w \in \mathbb{E}^n$.   The smocking pseudometric $\overline{d} : \mathbb{E}^N \times \mathbb{E}^N \rightarrow [ 0, \infty )$ is defined as 

\[  \overline{d} (v,w)  = d ( \pi (v), \pi (w)) = \text{min} \{ d_k (v', w') : \pi (v) = \pi (v'), \pi (w) = \pi (w'), k \in \mathbb{N}  \}  .\]
    
\end{definition}

\noindent
When studying specific examples of smocked spaces one also refers to a set of quantities called the smocking constants.  We will not need these in this work so we refer the reader to [1] for details, but note that one also assumes a positive smocking separation factor when defining a smocked space.

We will also require the definition of a stable norm.  Let $G$ be a connected weighted graph equipped with a free cocompact action of $\mathbb Z^n$ and let $o$ be a chosen base point of the graph $G$. 

\begin{definition}
The stable norm is defined by

\[
\|v\|_{\mathrm{st}}
=
\lim_{k\to\infty}
\frac{1}{k}
d_G(o,kv\cdot o).
\]
\end{definition}

\noindent
The significance of these norms is that they provide a way to pass from the local geometry of a periodic weighted graph $G$ to its asymptotic large-scale structure.  Intuitively, the stable norm measures the average distance required to travel in each direction $v \in \mathbb{Z}^n$.  Under the action of $\mathbb{Z}^n$, this asymptotic metric induces a norm on $\mathbb{R}^n$, converting the geometry of the graph into that of a normed vector space which captures the original large-scale geometry.

We are interested here in working with a construction where the definition of a smocking is adapted to such graphs.  This motivates the following definition:

\begin{definition}
A graph-adapted smocked realisation of $G$ consists of a smocked space
$(X,d_X)$ together with a surjective map
\[
\pi:X\to G
\]
satisfying the following properties:

\begin{enumerate}

\item[(A1)] There exists $D>0$ such that
\[
\operatorname{diam}\bigl(\pi^{-1}(p)\bigr)\le D
\]
for every $p\in G$.

\item[(A2)] For all $x,y\in X$,
\[
d_G(\pi(x),\pi(y))
\le d_X(x,y)+D.
\]

\item[(A3)] Every geodesic segment $\gamma$ in $G$ joining
$p,q\in G$ admits a lift joining points
$x\in\pi^{-1}(p)$ and $y\in\pi^{-1}(q)$
whose length is at most $L(\gamma)+D$.
\end{enumerate}
\end{definition}

\noindent
Intuitively, the role of the graph $G$ is to encode the large-scale combinatorial structure created by the stitches on the smocked space.  The map $\pi$ can therefore be thought of as collapsing the local geometry of $X$ onto a periodic graph.  Conditions (A1)-(A3) express the fact that this collapse process only changes distances by a uniformly bounded amount: each graph vertex is represented by a uniformly bounded region of the smocked space and graph geodesics can be lifted with uniformly bounded additional length.

\section{Proof of Theorem 1}

\noindent
The first ingredient is a realization theorem from the theory of stable
norms on periodic graphs. Recall that a norm $F$ on $\mathbb R^n$ is
called rational polyhedral if its unit ball is a centrally symmetric
polytope whose vertices have rational coordinates.  
 Realization results for rational polyhedral norm balls appear
in work of Babenko-Balacheff and Jotz, and are closely related to the theory of stable norms
and homogenization for periodic metrics and periodic graphs developed by
Burago-Ivanov-Kleiner and Kotani-Sunada [4-7]. For the purposes of this paper we assume the
following realization statement.

\begin{lemma} \label{thm:stable_realisation}
Let $F$ be a centrally symmetric rational polyhedral norm on
$\mathbb R^n$. Then there exists a connected $\mathbb Z^n$--periodic
weighted graph $(G,d_G)$ whose stable norm is exactly $F$.
\end{lemma}

\noindent
We must next transfer the
large-scale geometry of $G$ to an associated smocked space.  If we assume a graph-adapted smocked realisation as in Definition 4, we may prove the following lemma:

\begin{lemma}
\label{prop:additive}
Let $(X,d_X)$ be a graph-adapted smocked realisation of $G$.  Then there exists a constant $C>0$ such that
\[
\bigl|
d_X(x,y)-d_G(\pi(x),\pi(y))
\bigr|
\le C
\]
for all $x,y\in X$.
\end{lemma}

\begin{proof}
Let $x,y\in X$ and write
\[
p=\pi(x),\qquad q=\pi(y).
\]

\noindent
We first estimate $d_G(p,q)$ in terms of $d_X(x,y)$.  By property {\rm(A2)},
\[
d_G(p,q)\le d_X(x,y)+D
\]
and hence
\[
d_G(\pi(x),\pi(y))-d_X(x,y)\le D.
\]

\noindent
For the reverse inequality, let
\[
\gamma:[0,1]\to G
\]
be a geodesic segment joining $p$ to $q$. Then
\[
L(\gamma)=d_G(p,q).
\]
By property {\rm(A3)}, there exist points
\[
x'\in \pi^{-1}(p), \qquad
y'\in \pi^{-1}(q)
\]
and a path $\widetilde{\gamma}$ in $X$ joining $x'$ to $y'$ such that
\[
L(\widetilde{\gamma})
\le d_G(p,q)+D.
\]
Since $x$ and $x'$ lie in the same fibre of $\pi$, property {\rm(A1)}
implies
\[
d_X(x,x')\le D.
\]
Similarly,
\[
d_X(y,y')\le D.
\]
Concatenating a path from $x$ to $x'$, the lifted path
$\widetilde{\gamma}$, and a path from $y'$ to $y$, we obtain a path
joining $x$ to $y$ of length at most
\[
D+(d_G(p,q)+D)+D.
\]
Therefore
\[
d_X(x,y)\le d_G(p,q)+3D,
\]
and hence
\[
d_X(x,y)-d_G(\pi(x),\pi(y))
\le 3D.
\]
Combining the two inequalities yields
\[
-D
\le
d_X(x,y)-d_G(\pi(x),\pi(y))
\le
3D.
\]
As a result, we have
\[
\bigl|d_X(x,y)-d_G(\pi(x),\pi(y))\bigr|
\le 3D.
\]
Taking $C=3D$ completes the proof.
\end{proof}

\noindent
This lemma then enables us to prove our main result.

\begin{proof}[Proof of Theorem 1] 
By Lemma 2, there exists a constant $C>0$ such that
\[
|d_X(x,y)-d_G(\pi(x),\pi(y))|
\le C
\]
for all $x,y\in X$.  Consequently, after rescaling by $\lambda^{-1}$
\[
\left|
\lambda^{-1}d_X(x,y)-
\lambda^{-1}d_G(\pi(x),\pi(y))
\right|
\le \frac{C}{\lambda}.
\]
It follows that the pointed Gromov-Hausdorff distance between
\[
(X,\lambda^{-1}d_X,x_0)
\]
and
\[
(G,\lambda^{-1}d_G,\pi(x_0))
\]
tends to zero as $\lambda\to\infty$, hence $X$ and $G$ have the same asymptotic cone.  Since $G$ is a $\mathbb Z^n$-periodic weighted graph with stable norm
$F$, the classical stable norm convergence theorem for periodic metrics
implies that
\[
(G,\lambda^{-1}d_G,\pi(x_0))
\xrightarrow[\lambda\to\infty]{\mathrm{pGH}}
(\mathbb R^n,F).
\]
Therefore
\[
(X,\lambda^{-1}d_X,x_0)
\xrightarrow[\lambda\to\infty]{\mathrm{pGH}}
(\mathbb R^n,F).
\]
Since the limit is independent of the chosen sequence
$\lambda\to\infty$, the tangent cone at infinity is unique and equal to
$(\mathbb R^n,F)$.
\end{proof}

\noindent
Corollary 2 shows that a large class of normed spaces occurs as tangent cones at infinity of graph-adapted smocked spaces. More importantly, Theorem \ref{thm:main} identifies a mechanism by which asymptotic geometry is transferred from periodic graph models to smocked spaces. Consequently, results from stable norm theory and periodic homogenization may be applied directly to the study of the large-scale geometry of these smocked spaces.

\noindent

\section*{Acknowledgments} This research was partly conducted whilst the author was visiting the Okinawa Institute of Science and 
Technology (OIST) through the Theoretical Sciences Visiting Program (TSVP).

\section*{Funding Information} The author acknowledges support from a London Mathematical Society Early Career Research Travel Grant (ECR-2526-51).      

\section*{Conflict of Interest} The author declares that he has no conflict of interest.

\section*{Ethical approval} Not applicable.

\section*{References}

\noindent
[1] C. Sormani et al., Smocked Metric Spaces and their Tangent Cones, Miss. J. Math. Sci. \textbf{33} (1), 27-99 (2021). 

\noindent
[2] V. Antonetti, M. Fahrazad and A. Yamin, The Checkered Smocked Space and its Tangent Cone.  arXiv:1912.06294.

\noindent
[3] C. Sormani et al., Conjectures on Convergence and Scalar Curvature, in M. Gromov and H.B. Lawson (ed.), \textit{Perspectives in Scalar Curvature} (World Scientific, 2023).

\noindent
[4] I. Babenko and F. Balacheff, Sur la forme de la boule unité de la norme stable unidimensionnelle, Manuscripta Math. 119, 347–358 (2006).

\noindent
[5] M. Jotz, Hedlund metrics and the stable norm, Diff. Geo. App. 27(4), 543-550 (2009).

\noindent
[6] D. Burago, S. Ivanov and B. Kleiner, On the structure of the stable norm of periodic metrics, Math. Res. Lett. 4(6), 791-808 (1997).

\noindent
[7] M. Kotani and T. Sunada, Albanese Maps and Off Diagonal Long Time Asymptotics for the Heat Kernel, Commun. Math. Phys. 209, 633–670 (2000).

\noindent

\noindent

\noindent

\noindent






\end{document}